\documentclass[a4paper,12pt]{amsart} 
 \usepackage{epsf}
 \usepackage{amscd}                  
 \usepackage{xypic}                  
 \usepackage{amssymb}
\begin{document}
\def\D{{\mathcal D}}
\def\A{{\mathcal A}}
\def\C{{\mathcal C}}

\def\E{{\mathcal E}}
\def\Q{{\mathcal Q}}
\def\P{{\mathcal P}}
\def\curlyl{{\prec}}
\def\curlyr{{\succ}}
\newtheorem{thm}{Theorem}[section]
\newtheorem{cor}[thm]{Corollary}
\newtheorem{conjecture}[thm]{Conjecture}
\newtheorem{lemma}[thm]{Lemma}
\newtheorem{prop}[thm]{Proposition}
\newtheorem{defi}[thm]{Definition}
\newtheorem{nota}[thm]{Notation}
\newtheorem{rem}[thm]{Remark}
\newtheorem{example}[thm]{Example}

\title {The operad Quad is Koszul}

\author { Jon Eivind Vatne }

\address{ Matematisk institutt\\
          Johs. Brunsgt. 12 \\
          N-5008 Bergen \\
          Norway}   
        
\email{ jonev@mi.uib.no}

\begin{abstract} 
The purpose of this paper is to prove the koszulity of the operad
${\mathcal Quad}$, governing quadri-algebras.  That 
${\mathcal Quad}$ is Koszul was conjectured by Aguiar
and Loday in \cite{AL03}, where it was introduced.  The operad
${\mathcal Dend}$, governing dendriform algebras, is known to be
Koszul, \cite{Lod01}, and ${\mathcal Quad}$ is its second black square
power.  We find a new complex, based on the associahedron, which captures the
structure of ${\mathcal Dend}$.  This complex behaves well with
respect to the black squaring process, and allows us to conclude.  Also,
this proves koszulity of higher black square powers of ${\mathcal Dend}$.
\end{abstract}

\maketitle

\section{Introduction}

Given two quadratic binary operads $\P$ and $\Q$, there is a new
quadratic binary operad $\P\mbox{{\tiny $\blacksquare$}}\Q$ defined by
taking pairs of
operations from $\P$ and $\Q$, and imposing pairs of relations.  For
properties of this operation, see Ebrahimi-Fard and Guo, \cite{KFL04}.
Consider for example the operad ${\mathcal Dend}$, governing
dendriform algebras.  It was defined by Loday in \cite{Lod01}.  We can
form 

\[{\mathcal Quad}={\mathcal Dend}\,\mbox{{\tiny $\blacksquare$}}\, {\mathcal Dend}.\]

This operad governs quadri-algebras, and was introduced by Aguiar and
Loday in \cite{AL03}.
The operation $\mbox{{\tiny $\blacksquare$}}$ is mimicked on Manin's black dot operation on
quadratic algebras.  For operads, the question of how this operation
relates to koszulness, is far from being understood.  We will show
that in the present case, koszulity of ${\mathcal Quad}$ can be
deduced from the koszulity of ${\mathcal Dend}$.  We work over an
algebraically closed field $k$ of characteristic zero throughout this
paper.\\

We begin by formulating the koszulity condition in a useful form.  We
use the operadic bar construction for this.\\

Then we investigate ${\mathcal Dend}$ more closely.  It has two binary
operations $\curlyl$ and $\curlyr$, and three quadratic relations.
This data shows that the operation $*=\curlyl+\curlyr$ is associative.
So we have an associative operation that splits in two, and the
associativity axiom splits in three.  On the next level, we find that
the associahedron ${\mathcal A}_4$ splits into four parts.  Similarly,
the asscoiahedron ${\mathcal A}_n$ splits into $n$ parts.  We then see
that the sum of $n$ copies of the chain complex of the associahedron is
a direct summand of the operadic (dual) bar complex for ${\mathcal Dend}$.
Both these complexes are acyclic; the chain complex since it comes from
a polytope, and the bar complex since ${\mathcal Dend}$ is known to be
Koszul.  We can choose homotopy equivalences between the
complexes.\\

Now the operad ${\mathcal Quad}$ has 4 binary operations, and nine
quadratic relations.  The sum of these four operations is an
associative operation; the associativity axiom is exactly the sum of
the nine quadratic relations.  It then turns out that the
associahedron ${\mathcal A}_4$ splits into sixteen parts, and in
general ${\mathcal A}_n$ splits into $n^2$ parts.  This proves that the
dimension of ${\mathcal Quad}^{\,!}(n)$ is at least $n^2$, and we know
the opposite inequality from \cite{AL03}.  Collecting things together,
we can now use pairs of homotopy equivalences from ${\mathcal Dend}$
to show that the sum of $n^2$ copies of the chain complex of the
associahedron is homotopy equivalent to the operadic bar complex of
${\mathcal Quad}$.  Since the chain complex is acyclic, the bar complex
is acyclic, and ${\mathcal Quad}$ is Koszul.  During the proof, we
also verify the numerical conjecture from \cite{AL03}.\\

Instead of ${\mathcal Quad}={\mathcal Dend}^{\,\mbox{{\tiny $\blacksquare$}} \,2}$, we can of
course also consider ${\mathcal Dend}^{\,\mbox{{\tiny $\blacksquare$}} \,m}$.  The case
$m=3$ has been studied by Leroux \cite{Ler}, under the name of
octo-algebras.  The same proof as
for $m=2$ shows that this operad, for all $m$, is Koszul, modulo a
generalization of a lemma from \cite{AL03}.  We prove this lemma in
the last section.\\

Since all operads in this paper come from non-symmetric operads, the
symmetric group action will be suppressed throughout.  To get the true
operads from what is written here, tensor each algebraic construction
by $k[\Sigma_n]$ in degree $n$. 

\section{Koszulity for operads}

The operads we consider in this paper are of the following special
form:  they are generated by a finite number of binary operations.
Their relations are quadratic, and take the form

\[(x\circ_1 y)\circ_2 z= x\circ_3(y\circ_4 z)\]

where $\circ_i$ are binary operations.  It will be convenient to think
of such a relation as a directed edge between labelled trees

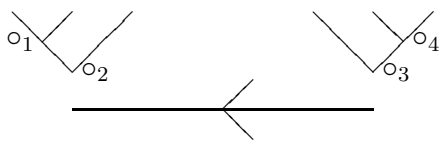
\[\begin{xy}
<0mm,0mm>;<40mm,0mm>**@{-},
<0mm,5mm>;<8mm,13mm>**@{-},
<0mm,5mm>;<-8mm,13mm>**@{-},
<-4mm,9mm>;<0mm,13mm>**@{-},
<40mm,5mm>;<48mm,13mm>**@{-},
<40mm,5mm>;<32mm,13mm>**@{-},
<44mm,9mm>;<40mm,13mm>**@{-},
<3mm,5mm>*{\circ_2},
<-7mm,9mm>*{\circ_1},
<43mm,5mm>*{\circ_3},
<47mm,9mm>*{\circ_4},
<20mm,0mm>;<24mm,4mm>**@{-},
<20mm,0mm>;<24mm,-4mm>**@{-},
\end{xy}\]

\begin{rem}
All trees considered here will have a finite number of leaves at the
top (3 in the two trees above), and some vertices below the leaves (2
in the examples above).  Each vertex will have a number ($\geq 2$) of
incoming edges.  The root is the lowest vertex
(the roots are labelled by $\circ_2$ and $\circ_3$ above).  Each
vertex apart from the root has a unique outoing edge.  The trees will
have labels at the vertices; these correspond to operations.  Labels
at the leaves correspond to inputs.  Since all the relations we will
consider have all the inputs in the same order, the labels at the
leaves will be suppressed throughout.  An edge between vertices will
be referred to as an {\em internal edge}, to distinguish it from a
leaf.
\end{rem}

We dub the space of binary operations $\Omega_{\P}$, and the space of
relations $\Lambda_{\P}$, following the notation from \cite{KFL04}. Note
that $\Lambda_{\P}\subset \Omega_{\P}^{\otimes 2}\oplus
\Omega_{\P}^{\otimes 2}$.  We write $\P=\P(\Omega_{\P},\Lambda_{\P})$.\\

Now we can define the squaring operation, still following \cite{KFL04}
(this has also been considered in \cite{Lod04}): 

\begin{defi}
The {\em black square product} of two operads
$\P=\P(\Omega_{\P},\Lambda_{\P})$ and $\Q=\P(\Omega_{\Q},\Lambda_{\Q})$
is

\[\P\,\mbox{{\tiny $\blacksquare$}}\,\Q=\P(\Omega_{\P}\otimes \Omega_{\Q},S_{23}(\Lambda_{\P}\otimes
\Lambda_{\Q})).\]

The operator $S_{23}$ simply switches tensor factors, so that the relations
come at the right place.
\end{defi}

We also need the quadratic dual algebra
$\P^!=\P(\Omega_{\P}^{\vee},\Lambda_{\P}^{\perp})$.  Here
$\Omega_{\P}^{\vee}$ is the linear dual tensored with the sign
representation (since we suppress the symmetric group action, this
merely involves a sign in the pairing), and the perpendicular is with
respect to a natural pairing; see {\em e.g.} Loday \cite{Lod96} for details.\\

To an operad we can associate its bar complex; this is basically the
free operad on the linear dual of the operad.  It is a dg operad; for
a quadratic operad, its zeroeth homology is the quadratic dual.  See
Markl, Snhider and Stasheff \cite{MSS} for details, including the
grading convention.  Since it will be important for us, we will
give the explicit structure of this construction for the operad
${\mathcal Dend}$ (actually its dual) later.  As part of the proof of
the main theorem, we will also find in explicit form the (dual) bar
complex of the operad ${\mathcal Quad}$.\\

\begin{defi}\label{koszul}
An operad (quadratic, binary) is called {\em Koszul} if the bar
complex is a resolution of the dual operad.
\end{defi}

\section{The higher degree structure of ${\mathcal Dend}$; splitting
  the associahedron}

We begin by writing out the dual bar complex of the operad ${\mathcal
  Dend}$ ({\em i.e.} the bar complex of the dual operad ${\mathcal
  Dias}={\mathcal Dend}^{\, !}$).  This operad, introduced in
\cite{Lod01}, governs dialgebras.  Then we make a complex out of the
associahedron, and finally we link these two together.\\

\subsection{The dual bar complex of ${\mathcal Dend}$}

\begin{defi}
The operad ${\mathcal Dend}$ is generated by two binary operations $\curlyl$
and $\curlyr$ satisfying three axioms.  In the language of trees, the
axioms can be written (let $*=\curlyl+\curlyr$)

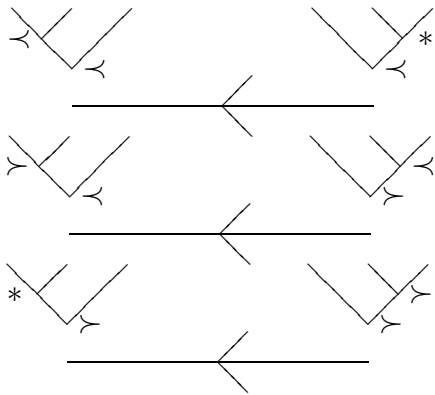
\[\begin{array}{c}
\begin{xy}
<0mm,0mm>;<40mm,0mm>**@{-},
<0mm,5mm>;<8mm,13mm>**@{-},
<0mm,5mm>;<-8mm,13mm>**@{-},
<-4mm,9mm>;<0mm,13mm>**@{-},
<40mm,5mm>;<48mm,13mm>**@{-},
<40mm,5mm>;<32mm,13mm>**@{-},
<44mm,9mm>;<40mm,13mm>**@{-},
<3mm,5mm>*{\curlyl},
<-7mm,9mm>*{\curlyl},
<43mm,5mm>*{\curlyl},
<47mm,9mm>*{*},
<20mm,0mm>;<24mm,4mm>**@{-},
<20mm,0mm>;<24mm,-4mm>**@{-},
\end{xy}\\

\begin{xy}
<0mm,0mm>;<40mm,0mm>**@{-},
<0mm,5mm>;<8mm,13mm>**@{-},
<0mm,5mm>;<-8mm,13mm>**@{-},
<-4mm,9mm>;<0mm,13mm>**@{-},
<40mm,5mm>;<48mm,13mm>**@{-},
<40mm,5mm>;<32mm,13mm>**@{-},
<44mm,9mm>;<40mm,13mm>**@{-},
<3mm,5mm>*{\curlyl},
<-7mm,9mm>*{\curlyr},
<43mm,5mm>*{\curlyr},
<47mm,9mm>*{\curlyl},
<20mm,0mm>;<24mm,4mm>**@{-},
<20mm,0mm>;<24mm,-4mm>**@{-},
\end{xy}\\

\begin{xy}
<0mm,0mm>;<40mm,0mm>**@{-},
<0mm,5mm>;<8mm,13mm>**@{-},
<0mm,5mm>;<-8mm,13mm>**@{-},
<-4mm,9mm>;<0mm,13mm>**@{-},
<40mm,5mm>;<48mm,13mm>**@{-},
<40mm,5mm>;<32mm,13mm>**@{-},
<44mm,9mm>;<40mm,13mm>**@{-},
<3mm,5mm>*{\curlyr},
<-7mm,9mm>*{*},
<43mm,5mm>*{\curlyr},
<47mm,9mm>*{\curlyr},
<20mm,0mm>;<24mm,4mm>**@{-},
<20mm,0mm>;<24mm,-4mm>**@{-},
\end{xy}
\end{array}\]

Note that the sum of these three relations is the associativity of $*$.

\end{defi}

\begin{lemma}\label{bardend}

\begin{itemize}
\item[(i)] The dual bar complex ${\tilde \D}=\D({\mathcal
    Dend}^!)=\D({\mathcal Dias})$ has the following graded parts:\\

\[\begin{array}{l}
\D_2^0\leftarrow 0\\
\D_3^0\leftarrow \D_3^{-1}\leftarrow 0\\
\D_4^0\leftarrow \D_4^{-1}\leftarrow \D_4^{-2}\leftarrow 0\\
\cdots\\
\D_n^0\leftarrow \cdots\leftarrow D_n^{-n+3}\leftarrow
\D_n^{-n+2}\leftarrow 0\\
\cdots \end{array}\]

The piece $D_i^j$ has basis given by labelled trees with $i$ leaves and $i+j-1$
vertices, with $l$ choices of labels for each vertex
with $l$ incoming edges.

\item[(ii)] The zeroeth homology of $\tilde{\D}$ is the operad
  ${\mathcal Dend}$, the higher homology vanishes.
\end{itemize}
\end{lemma}

\begin{proof}
By definition, $\tilde{\D}$ is the free operad construction on the twisted
linear dual of ${\mathcal Dias}$; this means that the piece $D_i^j$ is
given by trees with $i$ leaves and $i+j-1$ vertices, where
each vertex with $l$ incoming edges is labelled by an element of a
vector space of the same dimension as ${\mathcal Dias}(l)$, see
\cite{MSS}.  This space has dimension $l$ (see \cite{Lod01}).\\

The differential of such a labelled tree $T$ can be understood
inductively.  First, an unlabelled tree $T'$ with one vertex less than
$T$ has a labelling which appears with non-zero coefficient in the
differential of $T$ if and only if $T$ is the result of contracting
an internal edge of $T'$.  In this case, there is a unique labelling with this
property.  For all vertices except the two vertices of the contracted
edge, the labelling is unchanged.  Say that the vertex of $T$ that is
the image of the contracted edge has $l$ incoming edges.  Then the
labelling is induced from the map $\D_{l}^{-l+3}\leftarrow
D_l^{-l+2}$.  The description of this map depends on an explicit
description of the basis, and will be given during the proof of
Proposition \ref{dendass}.\\

The second part is the definition of koszulity as in Definition
\ref{koszul}.  ${\mathcal Dend}$ (and thus ${\mathcal Dias}$) is
Koszul by \cite{Lod01}. \\
\end{proof}

\begin{defi}
The augmented dual bar complex of ${\mathcal Dend}$ is the complex
$\D$ which is equal to ${\tilde \D}$ in non-positive degrees, but is
augmented by ${\mathcal Dend}$ in degree $+1$.  It is thus exact
everywhere.
\end{defi}

\subsection{The associahedron}

Let $\A_n$ be the associahedron for $n$ inputs.  This is a polytopal
cell complex of
dimension $n-2$; in particular, its (augmented) chain complex
${\mathcal CA}_n$ is exact. We grade it by minus the dimension of the
cells, so it has the explicit form

\[\begin{array}{ll}
n=2: &  {\mathcal CA}^{1}_2\leftarrow {\mathcal CA}^0_2\leftarrow 0\\
n=3: & {\mathcal CA}^{1}_3\leftarrow {\mathcal CA}^0_3\leftarrow
{\mathcal CA}^{-1}_3\leftarrow 0\\
& \cdots\\
n: & {\mathcal CA}^{1}_n\leftarrow {\mathcal CA}^0_n\leftarrow \cdots
\leftarrow
{\mathcal CA}^{-n+2}_n\leftarrow 0\end{array}\]

A basis for ${\mathcal CA}^j_i$ is given by trees with $j$ leaves and
$i+j-1$ vertices (except for the case $j=1$, where the generator is
represented by the empty cell).  The differential of a tree $T$ has nonzero
coefficient in a tree $T'$, with one vertex less, if and only if $T$
is the result of contracting an internal edge of $T'$.  The
coefficient is then $\pm 1$, depending on the choice of
orientation of the associahedron.

\begin{rem}
${\mathcal CA}$ is the bar complex of the operad ${\mathcal Ass}$ governing
associative algebras.  Its exactness is equivalent to the koszulity of
${\mathcal Ass}$.
\end{rem}

Later on, it will be convenient for us to label each vertex of each
tree in the chosen basis for this chain complex by $*$.\\

\begin{example}[$\A_4$]

The associahedron for four inputs is a pentagon:

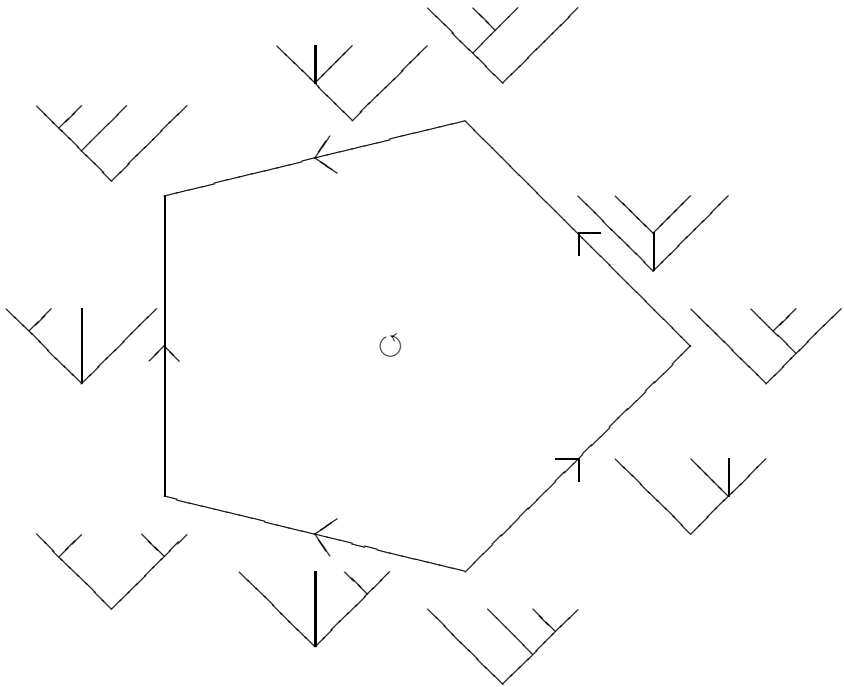
\[\begin{xy}
<30mm,20mm>*{\circlearrowleft},
<0mm,0mm>;<40mm,-10mm>**@{-},
<20mm,-5mm>;<23mm,-3mm>**@{-},
<20mm,-5mm>;<22mm,-8mm>**@{-},
<40mm,-10mm>;<70mm,20mm>**@{-},
<55mm,5mm>;<52mm,5mm>**@{-},
<55mm,5mm>;<55mm,2mm>**@{-},
<70mm,20mm>;<40mm,50mm>**@{-},
<55mm,35mm>;<58mm,35mm>**@{-},
<55mm,35mm>;<55mm,32mm>**@{-},
<40mm,50mm>;<0mm,40mm>**@{-},
<20mm,45mm>;<22mm,48mm>**@{-},
<20mm,45mm>;<23mm,43mm>**@{-},
<0mm,40mm>;<0mm,0mm>**@{-},
<0mm,20mm>;<2mm,18mm>**@{-},
<0mm,20mm>;<-2mm,18mm>**@{-},
<-7mm,42mm>;<-17mm,52mm>**@{-},
<-7mm,42mm>;<3mm,52mm>**@{-},
<-11mm,46mm>;<-5mm,52mm>**@{-},
<-14mm,49mm>;<-11mm,52mm>**@{-},
<-11mm,15mm>;<-21mm,25mm>**@{-},
<-11mm,15mm>;<-1mm,25mm>**@{-},
<-11mm,15mm>;<-11mm,25mm>**@{-},
<-18mm,22mm>;<-15mm,25mm>**@{-},
<-7mm,-15mm>;<-17mm,-5mm>**@{-},
<-7mm,-15mm>;<3mm,-5mm>**@{-},
<-14mm,-8mm>;<-11mm,-5mm>**@{-},
<0mm,-8mm>;<-3mm,-5mm>**@{-},
<20mm,-20mm>;<10mm,-10mm>**@{-},
<20mm,-20mm>;<30mm,-10mm>**@{-},
<20mm,-20mm>;<20mm,-10mm>**@{-},
<27mm,-13mm>;<24mm,-10mm>**@{-},
<45mm,-25mm>;<35mm,-15mm>**@{-},
<45mm,-25mm>;<55mm,-15mm>**@{-},
<49mm,-21mm>;<43mm,-15mm>**@{-},
<52mm,-18mm>;<49mm,-15mm>**@{-},
<70mm,-5mm>;<60mm,5mm>**@{-},
<70mm,-5mm>;<80mm,5mm>**@{-},
<75mm,0mm>;<70mm,5mm>**@{-},
<75mm,0mm>;<75mm,5mm>**@{-},
<80mm,15mm>;<70mm,25mm>**@{-},
<80mm,15mm>;<90mm,25mm>**@{-},
<84mm,19mm>;<78mm,25mm>**@{-},
<81mm,22mm>;<84mm,25mm>**@{-},
<65mm,30mm>;<55mm,40mm>**@{-},
<65mm,30mm>;<75mm,40mm>**@{-},
<65mm,30mm>:<65mm,35mm>**@{-},
<65mm,35mm>;<60mm,40mm>**@{-},
<65mm,35mm>;<70mm,40mm>**@{-},
<45mm,55mm>;<35mm,65mm>**@{-},
<45mm,55mm>;<55mm,65mm>**@{-},
<41mm,59mm>;<47mm,65mm>**@{-},
<44mm,62mm>;<41mm,65mm>**@{-},
<25mm,50mm>;<15mm,60mm>**@{-},
<25mm,50mm>;<35mm,60mm>**@{-},
<20mm,55mm>;<20mm,60mm>**@{-},
<20mm,55mm>;<25mm,60mm>**@{-},
\end{xy}\]

If we use the orientation as shown, and take ordered bases for the
chain complex by starting in the upper left corner and going
counter-clockwise, we get the following explicit complex ${\mathcal
  CA}_4$:

\[\xymatrix{
k & & k^5\ar[ll]_{(1\,1\,1\,1\,1)} & & & & &
k^5 \ar[lllll]_{\left(\begin{array}{ccccc}
1  & 0  & 0  & 0  & 1 \\
-1 & 1  & 0  & 0  & 0 \\
0  & -1 & -1 & 0  & 0 \\
0  & 0  & 1  & -1 & 0 \\
0  & 0  & 0  & 1  & -1\end{array}\right)}&
k \ar[l]_{\left(\begin{array}{c}-1\\-1\\1\\1\\1\end{array}\right)}}\]

\end{example}

\subsection{Splitting the associahedron}
The asociahedron for two inputs is just a point, labelled with $*$.
We will regard the two operations of ${\mathcal Dend}$ as a splitting
of this associahedron into two parts:

\[\begin{xy}
<0mm,0mm>;<-5mm,5mm>**@{-},
<0mm,0mm>;<5mm,5mm>**@{-},
<2mm,0mm>*{*},
<6mm,2mm>*{=},
<14mm,0mm>;<9mm,5mm>**@{-},
<14mm,0mm>;<19mm,5mm>**@{-},
<16mm,0mm>*{\curlyl},
<21mm,2mm>*{+},
<28mm,0mm>;<23mm,5mm>**@{-},
<28mm,0mm>;<33mm,5mm>**@{-},
<31mm,0mm>*{\curlyr},
\end{xy}\]

Similarly, the associahedron $\A_3$ is

\[\begin{xy}
<0mm,0mm>;<40mm,0mm>**@{-},
<0mm,5mm>;<8mm,13mm>**@{-},
<0mm,5mm>;<-8mm,13mm>**@{-},
<-4mm,9mm>;<0mm,13mm>**@{-},
<40mm,5mm>;<48mm,13mm>**@{-},
<40mm,5mm>;<32mm,13mm>**@{-},
<44mm,9mm>;<40mm,13mm>**@{-},
<3mm,5mm>*{*},
<-7mm,9mm>*{*},
<43mm,5mm>*{*},
<47mm,9mm>*{*},
<20mm,0mm>;<24mm,4mm>**@{-},
<20mm,0mm>;<24mm,-4mm>**@{-},
\end{xy}\]

which splits as the sum of the three relations of ${\mathcal Dend}$.
We dub the three relations of ${\mathcal Dend}$

\[\begin{array}{c}\begin{xy}
<-22mm,4mm>;<-30mm,12mm>**@{-},
<-22mm,4mm>;<-22mm,12mm>**@{-},
<-22mm,4mm>;<-14mm,12mm>**@{-},
<-20mm,4mm>*{1},
<-13mm,6mm>*{=},
<0mm,0mm>;<40mm,0mm>**@{-},
<0mm,5mm>;<8mm,13mm>**@{-},
<0mm,5mm>;<-8mm,13mm>**@{-},
<-4mm,9mm>;<0mm,13mm>**@{-},
<40mm,5mm>;<48mm,13mm>**@{-},
<40mm,5mm>;<32mm,13mm>**@{-},
<44mm,9mm>;<40mm,13mm>**@{-},
<3mm,5mm>*{\curlyl},
<-7mm,9mm>*{\curlyl},
<43mm,5mm>*{\curlyl},
<47mm,9mm>*{*},
<20mm,0mm>;<24mm,4mm>**@{-},
<20mm,0mm>;<24mm,-4mm>**@{-},
\end{xy}\\
\begin{xy}
<-22mm,4mm>;<-30mm,12mm>**@{-},
<-22mm,4mm>;<-22mm,12mm>**@{-},
<-22mm,4mm>;<-14mm,12mm>**@{-},
<-20mm,4mm>*{2},
<-13mm,6mm>*{=},
<0mm,0mm>;<40mm,0mm>**@{-},
<0mm,5mm>;<8mm,13mm>**@{-},
<0mm,5mm>;<-8mm,13mm>**@{-},
<-4mm,9mm>;<0mm,13mm>**@{-},
<40mm,5mm>;<48mm,13mm>**@{-},
<40mm,5mm>;<32mm,13mm>**@{-},
<44mm,9mm>;<40mm,13mm>**@{-},
<3mm,5mm>*{\curlyl},
<-7mm,9mm>*{\curlyr},
<43mm,5mm>*{\curlyr},
<47mm,9mm>*{\curlyl},
<20mm,0mm>;<24mm,4mm>**@{-},
<20mm,0mm>;<24mm,-4mm>**@{-},
\end{xy}\\
\begin{xy}
<-22mm,4mm>;<-30mm,12mm>**@{-},
<-22mm,4mm>;<-22mm,12mm>**@{-},
<-22mm,4mm>;<-14mm,12mm>**@{-},
<-20mm,4mm>*{3},
<-13mm,6mm>*{=},
<0mm,0mm>;<40mm,0mm>**@{-},
<0mm,5mm>;<8mm,13mm>**@{-},
<0mm,5mm>;<-8mm,13mm>**@{-},
<-4mm,9mm>;<0mm,13mm>**@{-},
<40mm,5mm>;<48mm,13mm>**@{-},
<40mm,5mm>;<32mm,13mm>**@{-},
<44mm,9mm>;<40mm,13mm>**@{-},
<3mm,5mm>*{\curlyr},
<-7mm,9mm>*{*},
<43mm,5mm>*{\curlyr},
<47mm,9mm>*{\curlyr},
<20mm,0mm>;<24mm,4mm>**@{-},
<20mm,0mm>;<24mm,-4mm>**@{-},
\end{xy}

\end{array}\]

The fact that the sum of these three relations is the associativity
condition can be written as

\[\begin{xy}<-22mm,4mm>;<-30mm,12mm>**@{-},
<-22mm,4mm>;<-22mm,12mm>**@{-},
<-22mm,4mm>;<-14mm,12mm>**@{-},
<-20mm,4mm>*{1},
<-13mm,6mm>*{+},\end{xy}
\
\begin{xy}<-22mm,4mm>;<-30mm,12mm>**@{-},
<-22mm,4mm>;<-22mm,12mm>**@{-},
<-22mm,4mm>;<-14mm,12mm>**@{-},
<-20mm,4mm>*{2},
<-13mm,6mm>*{+},\end{xy}
\
\begin{xy}<-22mm,4mm>;<-30mm,12mm>**@{-},
<-22mm,4mm>;<-22mm,12mm>**@{-},
<-22mm,4mm>;<-14mm,12mm>**@{-},
<-20mm,4mm>*{3},
<-13mm,6mm>*{=},\end{xy}
\
\begin{xy}<-22mm,4mm>;<-30mm,12mm>**@{-},
<-22mm,4mm>;<-22mm,12mm>**@{-},
<-22mm,4mm>;<-14mm,12mm>**@{-},
<-20mm,4mm>*{*},
\end{xy}\]

The chain complex ${\mathcal CA}_3$ splits in the same way:

\[\xymatrix{
k^3\ar[d]^+ & k^6 \ar[l]\ar[d]^+ & k^3\ar[l]\ar[d]^+\\
k   & k^2 \ar[l]         & k\ar[l]
}
\]

\begin{prop}\label{dendahedron}
The chain complex ${\mathcal CA}_n$ splits as a sum of $n$ chain
complexes constructed from ${\mathcal Dend}$.  These can be labelled
from $1$ to $n$ by considering the labels of the tree

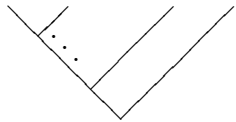
\[\begin{xy}
<0mm,0mm>;<-15mm,15mm>**@{-},
<0mm,0mm>;<15mm,15mm>**@{-},
<-4mm,4mm>;<7mm,15mm>**@{-},
<-11mm,11mm>;<-7mm,15mm>**@{-},
<-9mm,11mm>*{.},
<-7.5mm,9.5mm>*{.},
<-6mm,8mm>*{.},
\end{xy}\]

The first copy has labels $\curlyl,\cdots,\curlyl$ on this tree (from
top left), the second has $\curlyr,\curlyl,\cdots\curlyl$, the third
has $*,\curlyr,\curlyl,\cdots \curlyl$ and so on.  The two last have
labels $*,\cdots,*,\curlyr,\curlyl$ and $*,\cdots,*,\curlyr$.

\end{prop}

\begin{rem}
In particular, the binary operations are relabelled as

\[\begin{xy}
<0mm,0mm>;<-5mm,5mm>**@{-},
<0mm,0mm>;<5mm,5mm>**@{-},
<2mm,0mm>*{1},
<6mm,2mm>*{=},
<14mm,0mm>;<9mm,5mm>**@{-},
<14mm,0mm>;<19mm,5mm>**@{-},
<16mm,0mm>*{\curlyl},
\end{xy}
\
\, , \,
\
\begin{xy}
<0mm,0mm>;<-5mm,5mm>**@{-},
<0mm,0mm>;<5mm,5mm>**@{-},
<2mm,0mm>*{2},
<6mm,2mm>*{=},
<14mm,0mm>;<9mm,5mm>**@{-},
<14mm,0mm>;<19mm,5mm>**@{-},
<16mm,0mm>*{\curlyr},
\end{xy}\]
\end{rem}

\begin{proof}
We choose once and for all an orientation of each associahedron, and
an induced orientation on all cells, which
we will use consistently for each copy of it.  We do this so that all
edges are oriented as

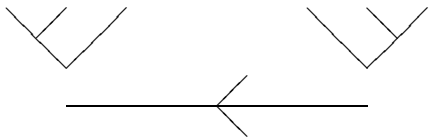
\[\begin{xy}
<0mm,0mm>;<40mm,0mm>**@{-},
<0mm,5mm>;<8mm,13mm>**@{-},
<0mm,5mm>;<-8mm,13mm>**@{-},
<-4mm,9mm>;<0mm,13mm>**@{-},
<40mm,5mm>;<48mm,13mm>**@{-},
<40mm,5mm>;<32mm,13mm>**@{-},
<44mm,9mm>;<40mm,13mm>**@{-},
<20mm,0mm>;<24mm,4mm>**@{-},
<20mm,0mm>;<24mm,-4mm>**@{-},
\end{xy}\]

Each edge is exactly represented by this move inside a larger tree.\\

  The proof starts with
describing what happens to the labels of each tree for the $n$ copies.
The sum over each tree gives the label of the associahedron.  Then we
do the same for the differentials.  This part is inductive.  Note
first that the proposition is consistent with what we have seen
already for $n=2,3$.  Thus the start of the induction is taken care of.\\

Let $T$ be a tree.  For each leaf, there is a unique path running
downwards from the leaf to the root.  The leftmost branch of the tree
is for instance the path from the leftmost leaf to the root.  In the
first copy, label all vertices along this leg by $1$, and all
remaining vertices by $*$.  Given two neighbouring leaves, say number
$i$ and $i+1$ from the left, there is
a unique shortest path connecting them.  We go from copy number $i$ to
copy number $i+1$ by changing the label of each vertex along this path
by the following rule:  If the label at a vertex is an integer $r$,
which is less than the number of incoming edges, increase it by one.
If the label is equal to the number of incoming edges, replace it by
$*$.  If the label is $*$, replace it by $1$.  In this way we finally end up
with a tree where all the vertices along the rightmost leg is labelled
by the number of incoming edges, all other vertices are labelled by
$*$.  For instance,

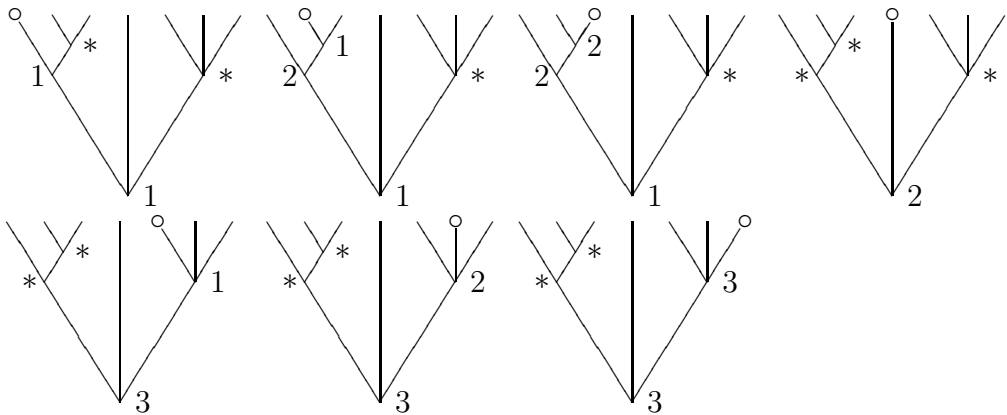
\[\begin{array}{llll}
\begin{xy}
<5mm,0mm>;<-5mm,16mm>**@{-},
<5mm,0mm>;<5mm,24mm>**@{-},
<5mm,0mm>;<15mm,16mm>**@{-},
<-5mm,16mm>;<-10mm,24mm>*{\circ}**@{-},
<-5mm,16mm>;<-2.5mm,20mm>**@{-},
<-2.5mm,20mm>;<-5mm,24mm>**@{-},
<-2.5mm,20mm>;<0mm,24mm>**@{-},
<15mm,16mm>;<10mm,24mm>**@{-},
<15mm,16mm>;<15mm,24mm>**@{-},
<15mm,16mm>;<20mm,24mm>**@{-},
<0mm,20mm>*{*},
<-7mm,16mm>*{1},
<8mm,0mm>*{1},
<18mm,16mm>*{*},
\end{xy}
&
\begin{xy}
<5mm,0mm>;<-5mm,16mm>**@{-},
<5mm,0mm>;<5mm,24mm>**@{-},
<5mm,0mm>;<15mm,16mm>**@{-},
<-5mm,16mm>;<-10mm,24mm>**@{-},
<-5mm,16mm>;<-2.5mm,20mm>**@{-},
<-2.5mm,20mm>;<-5mm,24mm>*{\circ}**@{-},
<-2.5mm,20mm>;<0mm,24mm>**@{-},
<15mm,16mm>;<10mm,24mm>**@{-},
<15mm,16mm>;<15mm,24mm>**@{-},
<15mm,16mm>;<20mm,24mm>**@{-},
<0mm,20mm>*{1},
<-7mm,16mm>*{2},
<8mm,0mm>*{1},
<18mm,16mm>*{*},
\end{xy}
&
\begin{xy}
<5mm,0mm>;<-5mm,16mm>**@{-},
<5mm,0mm>;<5mm,24mm>**@{-},
<5mm,0mm>;<15mm,16mm>**@{-},
<-5mm,16mm>;<-10mm,24mm>**@{-},
<-5mm,16mm>;<-2.5mm,20mm>**@{-},
<-2.5mm,20mm>;<-5mm,24mm>**@{-},
<-2.5mm,20mm>;<0mm,24mm>*{\circ}**@{-},
<15mm,16mm>;<10mm,24mm>**@{-},
<15mm,16mm>;<15mm,24mm>**@{-},
<15mm,16mm>;<20mm,24mm>**@{-},
<0mm,20mm>*{2},
<-7mm,16mm>*{2},
<8mm,0mm>*{1},
<18mm,16mm>*{*},
\end{xy}
&
\begin{xy}
<5mm,0mm>;<-5mm,16mm>**@{-},
<5mm,0mm>;<5mm,24mm>*{\circ}**@{-},
<5mm,0mm>;<15mm,16mm>**@{-},
<-5mm,16mm>;<-10mm,24mm>**@{-},
<-5mm,16mm>;<-2.5mm,20mm>**@{-},
<-2.5mm,20mm>;<-5mm,24mm>**@{-},
<-2.5mm,20mm>;<0mm,24mm>**@{-},
<15mm,16mm>;<10mm,24mm>**@{-},
<15mm,16mm>;<15mm,24mm>**@{-},
<15mm,16mm>;<20mm,24mm>**@{-},
<0mm,20mm>*{*},
<-7mm,16mm>*{*},
<8mm,0mm>*{2},
<18mm,16mm>*{*},
\end{xy}
\\
\begin{xy}
<5mm,0mm>;<-5mm,16mm>**@{-},
<5mm,0mm>;<5mm,24mm>**@{-},
<5mm,0mm>;<15mm,16mm>**@{-},
<-5mm,16mm>;<-10mm,24mm>**@{-},
<-5mm,16mm>;<-2.5mm,20mm>**@{-},
<-2.5mm,20mm>;<-5mm,24mm>**@{-},
<-2.5mm,20mm>;<0mm,24mm>**@{-},
<15mm,16mm>;<10mm,24mm>*{\circ}**@{-},
<15mm,16mm>;<15mm,24mm>**@{-},
<15mm,16mm>;<20mm,24mm>**@{-},
<0mm,20mm>*{*},
<-7mm,16mm>*{*},
<8mm,0mm>*{3},
<18mm,16mm>*{1},
\end{xy}
&
\begin{xy}
<5mm,0mm>;<-5mm,16mm>**@{-},
<5mm,0mm>;<5mm,24mm>**@{-},
<5mm,0mm>;<15mm,16mm>**@{-},
<-5mm,16mm>;<-10mm,24mm>**@{-},
<-5mm,16mm>;<-2.5mm,20mm>**@{-},
<-2.5mm,20mm>;<-5mm,24mm>**@{-},
<-2.5mm,20mm>;<0mm,24mm>**@{-},
<15mm,16mm>;<10mm,24mm>**@{-},
<15mm,16mm>;<15mm,24mm>*{\circ}**@{-},
<15mm,16mm>;<20mm,24mm>**@{-},
<0mm,20mm>*{*},
<-7mm,16mm>*{*},
<8mm,0mm>*{3},
<18mm,16mm>*{2},
\end{xy}
&
\begin{xy}
<5mm,0mm>;<-5mm,16mm>**@{-},
<5mm,0mm>;<5mm,24mm>**@{-},
<5mm,0mm>;<15mm,16mm>**@{-},
<-5mm,16mm>;<-10mm,24mm>**@{-},
<-5mm,16mm>;<-2.5mm,20mm>**@{-},
<-2.5mm,20mm>;<-5mm,24mm>**@{-},
<-2.5mm,20mm>;<0mm,24mm>**@{-},
<15mm,16mm>;<10mm,24mm>**@{-},
<15mm,16mm>;<15mm,24mm>**@{-},
<15mm,16mm>;<20mm,24mm>*{\circ}**@{-},
<0mm,20mm>*{*},
<-7mm,16mm>*{*},
<8mm,0mm>*{3},
<18mm,16mm>*{3},
\end{xy} & \end{array}\]

It is now obvious that the sum over all these labellings give the
labels $*$ at each spot, e.g. by induction (remove the root, and look
at the forest of smaller trees that remains).\\

As for the differential, we need to see what happens if the tree $T$ comes
from the tree $T'$ by contracting an internal edge, and the labels are as
prescribed at the $i$th level for both of them.  It is enough to
consider this in the case that $T$ has only one vertex, and $T'$ two.
If we forget about the labels, this means that $T$ represents the big
cell of the associahedron, whereas $T'$ represents a facet.  By our
choice of orientation, the differential of the tree $T$ has
coefficient $\pm 1$ on $T'$ in ${\mathcal CA}_n$.  In each copy of our
new complex, we use the same coefficient.  Then, when we sum over the
$n$ copies of the associahedron, we get that the differential of
${\mathcal CA}_n$ is the sum of the differential of the copies.  This
concludes the proof.

\end{proof}

\begin{defi}

We use the notation ${\mathcal DA}_n$ for the direct sum of the $n$
copies of the chain complex of the associahedron constructed in the
proposition.
\end{defi}

\begin{prop}\label{dendass}
There is a map from ${\mathcal DA}_n$ to the augmented dual bar
complex $\D_n$, which identifies ${\mathcal DA}_n$ with a direct
summand of $\D_n$. Each basis element of
${\tilde\D_n}$ appears with non-zero coefficient in the image of exactly one
basis element of ${\mathcal DA}_n$, where the coefficient is $1$.
\end{prop}

\begin{proof}
The map from ${\mathcal DA}_n$ to $\D_n$ is an isomorphism in degree
$-n+2$, where the two parts have the same dimension.  This gives us a
choice of basis for $\D_n$ as explained in the proof of Lemma
\ref{bardend}: a tree $T$ with labels on each vertex running from $1$
to the number of incoming edges.  This labelling corresponds to the
ordering from Proposition \ref{dendahedron}.  The map from ${\mathcal DA}_n$ to
$\D_n$ is given in general non-positive degrees by sending a labelled
tree $T$ to the tree
with the same labels in $\D$, understood as the sum where we split 
each label $*$ into the sum of the labels from $1$ to the number of
incoming edges of the vertex.  The differential of $T$ splits in the same
way; for each labelled tree $T'$ which appears with non-zero
coefficient in the differential of $T$, this coefficient is repeated
in the differential as many times as there are summands of $T'$.
 We extend the map to degree $1$ by taking the induced map on
cokernels.  Since all the coefficient of the inclusion map are $0$ or
$1$, we may choose a projection.  Now the statement of the proposition
is clear.
\end{proof}

\begin{rem}
Note that the description of the differential in $\D$, using the
explicit basis for $\D_n^{-n+2}$ given by the isomorphism with
${\mathcal DA}_n^{-n+2}$, fulfills the remaining part of the proof of
Lemma \ref{bardend}.
\end{rem}

\begin{defi}\label{htnot}

We define a few maps relating these two complexes:  first, write
$d_{\mathcal Dend}$ for the differential in $\D$.  Let the inclusion
of the summand be $f_{\mathcal Dend}$, the projection $p_{\mathcal
  Dend}$.  Then we let $h_{\mathcal Dend}$
be a homotopy between $p_{\mathcal Dend}:\D\rightarrow \D$ and the
identity map $I_{\mathcal Dend}$ on $\D$; this exists since it is projection
on a direct summand, and both complexes are split exact.  The homotopy
equivalence then takes the form

\[I_{\mathcal Dend}-p_{\mathcal Dend}=d_{\mathcal Dend}h_{\mathcal
  Dend}+h_{\mathcal Dend}d_{\mathcal Dend}\] 

\end{defi}

We will use these maps to construct similar maps for ${\mathcal Quad}$
in the next section.\\

\begin{rem}
The chain complex we have constructed is the chain complex of the
disjoint union of a number of copies of the associahedron, provided
that we include one ``empty cell'' for each copy.
\end{rem}

\section{Koszulity of ${\mathcal Quad}$}

Recall that

\[{\mathcal Quad}={\mathcal Dend}\,\mbox{{\tiny $\blacksquare$}} \,{\mathcal Dend}\]

by definition.  Using the column notation from \cite{KFL04}, this can
be written explicitly as follows:  There are four binary
operations  
\[\left[ \begin{array}{c} \curlyl \\ 
    \curlyl\end{array}\right],\,
\
 \left[ \begin{array}{c} \curlyl \\
    \curlyr\end{array}\right],\,
\
 \left[ \begin{array}{c} \curlyr \\
    \curlyl\end{array}\right] \,\mbox{and}\,
\ 
\left[ \begin{array}{c} \curlyr \\
    \curlyr\end{array}\right].\]

These satisfy nine relations, which are pairs of the relations from
${\mathcal Dend}$.  We label them

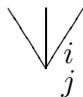
\[\begin{xy}
<0mm,0mm>;<5mm,8mm>**@{-},
<0mm,0mm>;<0mm,8mm>**@{-},
<0mm,0mm>;<-5mm,8mm>**@{-},
<3mm,2mm>*{i},
<3mm,-2mm>*{j},
\end{xy}\]

where $i$ and $j$ run from $1$ to $3$.  For instance, 

\[\begin{xy}
<0mm,0mm>;<5mm,8mm>**@{-},
<0mm,0mm>;<0mm,8mm>**@{-},
<0mm,0mm>;<-5mm,8mm>**@{-},
<3mm,2mm>*{1},
<3mm,-2mm>*{2},
\end{xy}
\
=
\
\begin{xy}
<0mm,0mm>;<40mm,0mm>**@{-},
<0mm,5mm>;<8mm,13mm>**@{-},
<0mm,5mm>;<-8mm,13mm>**@{-},
<-4mm,9mm>;<0mm,13mm>**@{-},
<40mm,5mm>;<48mm,13mm>**@{-},
<40mm,5mm>;<32mm,13mm>**@{-},
<44mm,9mm>;<40mm,13mm>**@{-},
<3mm,5mm>*{\curlyl},
<3mm,2mm>*{\curlyl},
<-7mm,9mm>*{\curlyl},
<-7mm,6mm>*{\curlyr},
<43mm,5mm>*{\curlyl},
<43mm,2mm>*{\curlyr},
<47mm,9mm>*{*},
<47mm,6mm>*{\curlyl},
<20mm,0mm>;<24mm,4mm>**@{-},
<20mm,0mm>;<24mm,-4mm>**@{-},
\end{xy}\]

This section is devoted to proving the main theorem of this paper:

\begin{thm}\label{Quadiskos}
The operad ${\mathcal Quad}$ is Koszul.
\end{thm}

The proof proceeds by mimicking the constructions we have made for
${\mathcal Dend}$ earlier, but reversing the final implication.

\begin{prop}\label{split}

The chain complex of the associahedron ${\mathcal CA}_n$ splits as a
sum of $n^2$ chain complexes constructed from ${\mathcal Quad}$.  These
can be labelled by indexes $i,j$ running from $1$ to $n$, where the
tree

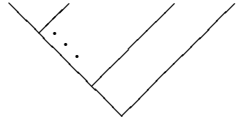
\[\begin{xy}
<0mm,0mm>;<-15mm,15mm>**@{-},
<0mm,0mm>;<15mm,15mm>**@{-},
<-4mm,4mm>;<7mm,15mm>**@{-},
<-11mm,11mm>;<-7mm,15mm>**@{-},
<-9mm,11mm>*{.},
<-7.5mm,9.5mm>*{.},
<-6mm,8mm>*{.},
\end{xy}\]

has pairs of labels as in Proposition \ref{dendahedron}.

\end{prop}

The proof is exactly as for Proposition \ref{dendahedron}.\\

\begin{defi}
The complex constructed from $n^2$ copies of the associahedron is
denoted by ${\mathcal QA}_n$; it is exact everywhere.
\end{defi}

\begin{defi}
The dual bar complex of ${\mathcal Quad}$ is denoted by ${\tilde \E}$, the
augmented version is denoted by $\E$.
\end{defi}

\begin{prop}\label{Dirsum}
\begin{itemize}
\item[(a)]There is a map from ${\mathcal QA}_n$ to $\E$, identifying
  ${\mathcal QA}_n$ with a direct summand.  Each basis
element of ${\tilde \E}$ appears with non-zero coefficient in the
image of exactly one basis element of ${\mathcal QA}$, where the
coefficient is $1$.
\item[(b)] The dimension of ${\mathcal Quad}^{\,!}(n)$ is $n^2$.
\end{itemize}
\end{prop}

The basis for ${\tilde \E}$ will be constructed in the proof.

\begin{proof} 
We will proceed by induction on $n$, the cases $n=2,3$ being trivial
(for both parts of the proposition).  Now by induction we have the
following diagram:

\[\xymatrix{
{\mathcal Quad}(n) & \E_n^0\ar[l] & \cdots\ar[l] &\E_n^{-n+3}\ar[l]
&\E_n^{-n+2}\ar[l] \\
{\mathcal QA}_n^1 & {\mathcal QA}_n^0\ar[l]\ar[u]&\cdots\ar[l] &{\mathcal
  QA}_n^{-n+3}\ar[l]\ar[u]& {\mathcal QA}_n^{-n+2}\ar[l]}\]

Each vertical arrow is the inclusion of a direct summand.  Now the
induced map from ${\mathcal QA}_n^{-n+2}$ to $\E_n^{-n+2}$ is clearly
injective.  In particular, $\dim \E_n^{-n+2}\geq n^2$.  The opposite
equality is Lemma \ref{Aguiar-Loday}.  Thus this map is an
isomorphism.  We use this isomorphism to choose basis for
$\E_n^{-n+2}$.  In particular, it is the inclusion of a direct
summand.  The induced map on the left hand side is a cokernel of an
inclusion of a direct summand; all in all, the complex ${\mathcal
  QA}_n$ is a direct summand of $\E_n$.\\
By the choice of basis, we see that for each tree $T$, the set of
labellings for $T$ giving basis elements of $\E_n$ is the tensor
square of the same for $\D_n$.  The analogous statement is obviously
true for ${\mathcal QA}$ and ${\mathcal DA}$, and the map from
${\mathcal QA}$ to $\E$ is locally, {\em i.e.} for each tree, the tensor square
of the corresponding map from ${\mathcal DA}$ to ${\D}$.  In
particular, each basis element of $\E$ appears with nonzero
coefficient in the image of a unique basis element in ${\mathcal QA}$,
and the coefficient is $1$.
\end{proof}

\begin{lemma}[Aguiar-Loday]\label{Aguiar-Loday}
The dimension of ${\mathcal Quad}^{\, !}(n)$ is $\leq n^2$.
\end{lemma}

This is taken from \cite{AL03}.

\begin{proof}[Proof of Theorem \ref{Quadiskos}]
Since we know
that the basis elements of ${\mathcal QA}$ and $\E$ are given by
pairs of basis elements of ${\mathcal DA}$ and $\D$, respectively, and
that the maps respect this, we can simply form pairs of homotopies as
well.  Explicitly, using the notation from Definiition \ref{htnot}, we
get that the map from ${\mathcal DA}$ to $\D$, and the projection onto
the summand, and the choice of homotopy are

\[f_{\mathcal Quad}=\left[\begin{array}{c}f_{\mathcal Dend}\\ f_{\mathcal
      Dend}\end{array}\right],\,\,
p_{\mathcal Quad}=\left[\begin{array}{c}p_{\mathcal Dend}\\ p_{\mathcal
      Dend}\end{array}\right],\,\,
h_{\mathcal Quad}=\left[\begin{array}{c}h_{\mathcal Dend}\\ h_{\mathcal
      Dend}\end{array}\right]\]

Since also 

\[d_{\mathcal Quad}=\left[\begin{array}{c}d_{\mathcal Dend}\\ d_{\mathcal
      Dend}\end{array}\right]\]

we get the homotopy relation

\[I_{\mathcal Quad}-p_{\mathcal Quad}=d_{\mathcal Quad}h_{\mathcal
  Quad}+h_{\mathcal Quad}d_{\mathcal Quad}\] 

from the corresponding relation for ${\mathcal Dend}$ used twice.
Thus the dual bar complex is homotopic to an exact complex, and is
therefore itself exact.  This proved that ${\mathcal Quad}$ is Koszul.

\end{proof}

\begin{cor}

The dimension of ${\mathcal Quad}(n)$ is $d_n$, where

\[d_n=\frac{1}{n}\sum_{j=n}^{j=2n-1}\binom{3n}{n+1+j}\binom{j-1}{j-n}\]

\end{cor}

This follows from the numerical data for the dual, and the koszulity;
see \cite{AL03}.

\section{Generalization}

If we consider ${\mathcal Dend}^{\,\mbox{{\tiny $\blacksquare$}} \,m}$ for general $m$, the
proof of the main theorem goes through modulo the generalization of
Lemma \ref{Aguiar-Loday}.  The aim of this section is to prove this
generalization.\\

\begin{thm}

For each $m\geq 1$, the operad

\[{\mathcal Dend}^{\,\mbox{{\tiny $\blacksquare$}}\, m}={\mathcal Dend}\,\mbox{{\tiny $\blacksquare$}} \cdots \mbox{{\tiny $\blacksquare$}}\,
{\mathcal Dend}\]

is Koszul.
\end{thm}

Operations in ${\mathcal Dend}^{\,\mbox{{\tiny $\blacksquare$}} \,m}$ are represented by
$m$-tuples of operations in ${\mathcal Dend}$, relations by $m$-tuples
of relations, and so forth.  The generalizations of Proposition
\ref{split} and Proposition \ref{Dirsum} go through with the same
proof, as does the concluding proof of the theorem; we only need to check the generalization of Lemma
\ref{Aguiar-Loday}.  We need a convention about the combinatorial
structures that will appear in the proof, summarized as

\begin{nota}

We consider an $m$-dimensional hypercube, which is composed of $3^m$
unit hypercubes.  Each of the consituent unit cubes has coordinates; an
$m$-tuple of elements from $\{1,2,3\}$.  The {\em corners} are the
unit cubes where no coordinate is equal to $2$.  $(3,3,\dots,3)$ is
the cube with highest coordinates.  For each subset of
$\{1,\dots,m\}$, say with $j$ elements, there is a subcube of
dimension $j$ where we only use the coordinates in the subset, and set
all other coordinates to $1$.

\end{nota}

\begin{lemma}

The quadratic dual operad $({\mathcal Dend}^{\,\mbox{{\tiny $\blacksquare$}}\, m})^!$ satisfies

\[\dim ({\mathcal Dend}^{\,\mbox{{\tiny $\blacksquare$}} \,m})^!(n)\leq n^m\]

\end{lemma}

\begin{proof}

The proof is also a direct generalization of Aguiar-Loday's proof of
Lemma \ref{Aguiar-Loday}, see \cite{AL03}.\\

We denote the operations in the dual operad by the same symbols as
we denote the operations in the original operad, that is as $m$-tuples
of linear combinations of $\curlyl$ and $\curlyr$.  So this is a space
of dimension $2^m$.  We choose representatives in degree three, one
for each relation in $ {\mathcal Dend}^{\,\mbox{{\tiny $\blacksquare$}} \,m}$, that is one for
each $m$-tuple of relations for ${\mathcal Dend}$.  These form a space
of dimension $3^m$.  We label the relations as column vectors, where
each element is $1,2$ or $3$ (as with the labelling in ${\mathcal
  Dend}$ from Section 3.3).  Then we choose representatives as
follows:  For each relation with no element equal to three, we use the
tree

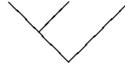
\[\begin{xy}
<0mm,5mm>;<8mm,13mm>**@{-},
<0mm,5mm>;<-8mm,13mm>**@{-},
<-4mm,9mm>;<0mm,13mm>**@{-},
\end{xy}\]

For each element with at least one element equal to three, we use the
tree

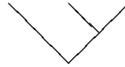
\[\begin{xy}
<0mm,5mm>;<8mm,13mm>**@{-},
<0mm,5mm>;<-8mm,13mm>**@{-},
<4mm,9mm>;<0mm,13mm>**@{-},
\end{xy}\]

For each element equal to $1$, we use the label $\curlyl$ at both
places, for each element equal to $3$ we use the label $\curlyr$ at
both places, and for each element equal to $2$ we use $\curlyr$ at the
leftmost vertex, $\curlyl$ at the rightmost vertex.  In particular,
whenever no element is equal to $2$, the upper and the lower label is
the same.  For example, with $m=3$ and the two vectors $(1,2,1)$ and
$(3,1,2)$ we get

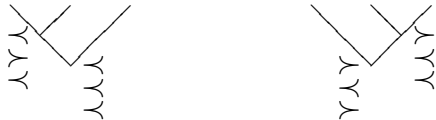
\[\begin{xy}
<0mm,5mm>;<8mm,13mm>**@{-},
<0mm,5mm>;<-8mm,13mm>**@{-},
<-4mm,9mm>;<0mm,13mm>**@{-},
<40mm,5mm>;<48mm,13mm>**@{-},
<40mm,5mm>;<32mm,13mm>**@{-},
<44mm,9mm>;<40mm,13mm>**@{-},
<3mm,5mm>*{\curlyl},
<3mm,2mm>*{\curlyl},
<3mm,-1mm>*{\curlyl},
<-7mm,9mm>*{\curlyl},
<-7mm,6mm>*{\curlyr},
<-7mm,3mm>*{\curlyl},
<37mm,5mm>*{\curlyr},
<37mm,2mm>*{\curlyl},
<37mm,-1mm>*{\curlyr},
<47mm,9mm>*{\curlyr},
<47mm,6mm>*{\curlyl},
<47mm,3mm>*{\curlyl},
\end{xy}\]

The set of elements where the local patterns (along each edge) is as
above, clearly generates the operad as a vector space.  Let $s_n$ be
the number of such elements.  We will show that $s_n=n^m$ by
recurrence, using these local patterns.  Specifically, we will write
$s_n$ as a sum of $2^m$ summands corresponding to the label at the
root; each of these summands can be written as a sum of terms with
lower degree, and this will give the recurrence.  Our explicit
knowledge about the situation in degree $2$ and $3$ gives the starting
point.\\

For each vector of labels, there is a unique corner of the $m$th
hypercube with that vector at the root.  This is true by the choice of
labellings.

We let $(\curlyl,\curlyl,\cdots,\curlyl)_n$ be the number of elements
of degree $n$ with the label
$(\curlyl,\cdots,\curlyl)$ at the root, $(\curlyr,\cdots,\curlyr)_n$
the number of elements with $(\curlyr,\cdots,\curlyr)$ at the root,
and a general $m$-vector of $\curlyl$s and $\curlyr$s, subscripted
$n$, represents the number of elements of degree $n$ with root
labelled by that vector.  Obviously, for any vector $(c_1,\cdots
,c_m)$ of such labels, we have $(c_1,\cdots,c_m)_2=1$.  This gives the
start of our recurrence.\\

The box with coordinates $(1,1,\cdots,1)$ has label
$(\curlyl,\curlyl,\cdots, \curlyl)$ at the root.  This is the same as
the label of the root for each box with no coordinate equal to $3$,
and each possible combination of labels appear exactly once at the
upper vertex of a tree in this hypercube (of size $2^m$).  So we get

\[(\curlyl,\cdots,\curlyl)_{n+1}=\sum_{c_i\in
  \{\curlyl,\curlyr\}}(c_1,\cdots,c_m)_n\]

In particular, $s_n=a_{n+1}$.\\

For each box with exactly one coordinate equal to $3$, all the rest
being one, the label of the root and the upper are equal, and no other
box has either of these labels at any vertex.  Thus

\[(c_1,\cdots,c_m)_n=1\]

if all $c_i$ are equal to $\curlyl$ except for one $\curlyr$.\\

In general, for each box with $j$ coordinates equal to $3$, the rest
being $1$, there is a subcube of dimension $j$, of size $3^j$, such
that the given box is the corner with highest coordinates in this
subcube.  Now
the label at the root of this box also appears on the root of all
the boxes in this subcube where all the coordinates are $2$ or $3$, with
at least one $3$.  So the recurrence relation for this box is equal to
the recurrence relation for vector with labels all $\curlyr$s in a hypercube
of dimension $j$.  We claim that this is

\[(\curlyr,\cdots,\curlyr)_n=\sum_{i=1}^j(-1)^{i-1}\binom{j}{i}(n-1)^{j-i}=
(n-1)^j-(n-2)^j\]

The case $j=1$ is the special case considered above, so we get the
start of the recurrence relation.\\

The recurrence relation in general says that

\[(\curlyr,\cdots,\curlyr)_{n+1}=\sum_{c_t\in\{\curlyl,\curlyr\},
  \exists \,t ,\,c_t\neq \curlyl}(c_1,\cdots,c_j)\]

By induction, the formula above holds for each summand on the right,
so we need to show that

\[n^j-(n-1)^j= \sum_{k=1}^j\binom{j}{k}((n-1)^k-(n-2)^k)\]

This follows from the Binomial theorem, on writing $n$ as $(n-1)+1$
and $n-1=(n-2)+1$; the two terms with $k=0$ cancel.\\

Note that there are $\binom{m}{j}$ vectors with $j$ $\curlyr$s, the
rest $\curlyl$s, so our final recurrence relation, for the vector
$(\curlyl,\cdots,\curlyl)$ takes the form 

\[(\curlyl,\cdots,\curlyl)_{n+1}=(\curlyl,\cdots,\curlyl)_n+\sum_{j=1}^{m}
\binom{m}{j}((n-1)^j-(n-2)^j)\]

This recurrence relation is satisfied by
$(\curlyl,\cdots,\curlyl)_{n+1}=n^m$; this also follows from the
Binomial theorem as above.\\

So we've shown that the vector space $({\mathcal Dend}^{\,\mbox{{\tiny
      $\blacksquare$}}\, m})^!(n)$ is spanned by $n^m$ elements; this
is therefore a bound on the dimension.

\end{proof}

\begin{rem}

The patterns we have chosen in the proof are modelled on the patterns
from \cite{AL03}, and for $m=2$ the proof reduces to their proof.  The
only difference is that we have chosen the other tree in position
$(3,1)$.

\end{rem}

\begin{rem}
Note that we have computed the dimension 

\[\dim {\mathcal Dend}^{\,\mbox{{\tiny $\blacksquare$}}\, m}(n)=n^m.\]

This follows from the proof of the theorem, exactly as in the case $m=2$.

\end{rem}

\end{document}